




%

%
 \font\twelvebf=cmbx12
 \font\twelvett=cmtt12
 \font\twelveit=cmti12
 \font\twelvesl=cmsl12
 \font\twelverm=cmr12		\font\ninerm=cmr9
 \font\twelvei=cmmi12		\font\ninei=cmmi9
 \font\twelvesy=cmsy10 at 12pt	\font\ninesy=cmsy9
 \skewchar\twelvei='177		\skewchar\ninei='177
 \skewchar\seveni='177	 	\skewchar\fivei='177
 \skewchar\twelvesy='60		\skewchar\ninesy='60
 \skewchar\sevensy='60		\skewchar\fivesy='60
%
%

%
 \font\fourteenrm=cmr12 scaled 1200
 \font\seventeenrm=cmr12 scaled 1440
 \font\fourteenbf=cmbx12 scaled 1200
 \font\seventeenbf=cmbx12 scaled 1440
%
%

%
%
%
\font\tenmsb=msbm10
\font\twelvemsb=msbm10 scaled 1200
\newfam\msbfam

%
\font\tensc=cmcsc10
\font\twelvesc=cmcsc10 scaled 1200
\newfam\scfam

%
\def\seventeenpt{\def\rm{\fam0\seventeenrm}%
 \textfont\bffam=\seventeenbf	\def\bf{\fam\bffam\seventeenbf}}
\def\fourteenpt{\def\rm{\fam0\fourteenrm}%
 \textfont\bffam=\fourteenbf	\def\bf{\fam\bffam\fourteenbf}}
\def\twelvept{\def\rm{\fam0\twelverm}%
 \textfont0=\twelverm	\scriptfont0=\ninerm	\scriptscriptfont0=\sevenrm
 \textfont1=\twelvei	\scriptfont1=\ninei	\scriptscriptfont1=\seveni
 \textfont2=\twelvesy	\scriptfont2=\ninesy	\scriptscriptfont2=\sevensy
 \textfont3=\tenex	\scriptfont3=\tenex	\scriptscriptfont3=\tenex
 \textfont\itfam=\twelveit	\def\it{\fam\itfam\twelveit}%
 \textfont\slfam=\twelvesl	\def\sl{\fam\slfam\twelvesl}%
 \textfont\ttfam=\twelvett	\def\tt{\fam\ttfam\twelvett}%
 \scriptfont\bffam=\tenbf 	\scriptscriptfont\bffam=\sevenbf
 \textfont\bffam=\twelvebf	\def\bf{\fam\bffam\twelvebf}%
 \textfont\scfam=\twelvesc	\def\sc{\fam\scfam\twelvesc}%
 \textfont\msbfam=\twelvemsb	
 \baselineskip 14pt%
 \abovedisplayskip 7pt plus 3pt minus 1pt%
 \belowdisplayskip 7pt plus 3pt minus 1pt%
 \abovedisplayshortskip 0pt plus 3pt%
 \belowdisplayshortskip 4pt plus 3pt minus 1pt%
 \parskip 3pt plus 1.5pt
 \setbox\strutbox=\hbox{\vrule height 10pt depth 4pt width 0pt}}
\def\tenpt{\def\rm{\fam0\tenrm}%
 \textfont0=\tenrm	\scriptfont0=\sevenrm	\scriptscriptfont0=\fiverm
 \textfont1=\teni	\scriptfont1=\seveni	\scriptscriptfont1=\fivei
 \textfont2=\tensy	\scriptfont2=\sevensy	\scriptscriptfont2=\fivesy
 \textfont3=\tenex	\scriptfont3=\tenex	\scriptscriptfont3=\tenex
 \textfont\itfam=\tenit		\def\it{\fam\itfam\tenit}%
 \textfont\slfam=\tensl		\def\sl{\fam\slfam\tensl}%
 \textfont\ttfam=\tentt		\def\tt{\fam\ttfam\tentt}%
 \scriptfont\bffam=\sevenbf 	\scriptscriptfont\bffam=\fivebf
 \textfont\bffam=\tenbf		\def\bf{\fam\bffam\tenbf}%
 \textfont\scfam=\tensc		\def\sc{\fam\scfam\tensc}%
 \textfont\msbfam=\tenmsb	
 \baselineskip 12pt%
 \abovedisplayskip 6pt plus 3pt minus 1pt%
 \belowdisplayskip 6pt plus 3pt minus 1pt%
 \abovedisplayshortskip 0pt plus 3pt%
 \belowdisplayshortskip 4pt plus 3pt minus 1pt%
 \parskip 2pt plus 1pt
 \setbox\strutbox=\hbox{\vrule height 8.5pt depth 3.5pt width 0pt}}

%
\def\twelvepoint{%
 \def\small{\tenpt\rm}%
 \def\normal{\twelvept\rm}%
 \def\large{\fourteenpt\rm}%
 \def\huge{\seventeenpt\rm}%
 \footline{\hss\twelverm\folio\hss}%
 \normal}
\def\tenpoint{%
 \def\small{\tenpt\rm}%
 \def\normal{\tenpt\rm}%
 \def\large{\twelvept\rm}%
 \def\huge{\fourteenpt\rm}%
 \footline{\hss\tenrm\folio\hss}%
 \normal}

\tenpoint

%

%
\catcode`\@=11
%
%
\def\footnote#1{\edef\@sf{\spacefactor\the\spacefactor}#1\@sf
 \insert\footins\bgroup\small
 \interlinepenalty100	\let\par=\endgraf
 \leftskip=0pt		\rightskip=0pt
 \splittopskip=10pt plus 1pt minus 1pt	\floatingpenalty=20000
 \smallskip\item{#1}\bgroup\strut\aftergroup\@foot\let\next}
%
%
%
%
\def\hexnumber@#1{\ifcase#1 0\or 1\or 2\or 3\or 4\or 5\or 6\or 7\or 8\or
 9\or A\or B\or C\or D\or E\or F\fi}
\edef\msbfam@{\hexnumber@\msbfam}
\def\Bbb#1{\fam\msbfam\relax#1}
%
%
%
\catcode`\@=12

\twelvepoint

\font\twelvebf=cmbx12	
\font\ninebf=cmbx9	\font\sevenbf=cmbx7	\font\fivebf=cmbx5

\font\twelvebfit=cmbxti10 at 12pt	
\newfam\bfitfam

\font\twelvebm=cmmib10 at 12pt		\font\tenbm=cmmib10
\font\ninebm=cmmib9	\font\sevenbm=cmmib7	\font\fivebm=cmmib5

\skewchar\twelvebm='177	\skewchar\tenbm='177
\skewchar\ninebm='177	\skewchar\sevenbm='177	\skewchar\fivebm='177

\def\twelvepointbold{\def\bold{
\textfont0=\twelvebf	\scriptfont0=\ninebf	\scriptscriptfont0=\sevenbf
\textfont1=\twelvebm	\scriptfont1=\ninebm	\scriptscriptfont1=\sevenbm
\textfont\bffam=\twelvebf	\textfont\bfitfam=\twelvebfit
\def\rm{\fam\bffam\twelvebf}%
\def\it{\fam\bfitfam\twelvebfit}%
\rm}}

\twelvepointbold
\font\twelvebsy=cmbsy10 scaled 1200
\font\fourteenbsy=cmbsy10 scaled 1440
\font\fourteenbf=cmbx12 scaled 1200
\def\Bold{\bold\textfont2=\twelvebsy}
\def\largebold{\textfont0=\fourteenbf\textfont2=\fourteenbsy}

\chardef\tempcat=\the\catcode`\@
\catcode`\@=11


\def\cydot{{\mathsurround=0pt$\cdot$}}

\def\ubar#1{\oalign{#1\crcr\hidewidth
    \vbox to.2ex{\hbox{\char22}\vss}\hidewidth}}

\def\cprime{\/{\mathsurround=0pt$'$}}
\def\Cprime{{\mathsurround=0pt$'$}}
\def\cdprime{\/{\mathsurround=0pt$''$}}
\def\Cdprime{{\mathsurround=0pt$\ubar{\hbox{$''$}}$}}

\def\dbar{dj}           
\def\Dbar{Dj}           

\def\dz{dz}
\def\Dz{Dz}
\def\dzh{dzh\cydot }
\def\Dzh{Dzh\cydot }


\def\@gobble#1{}
\def\@testgrave{\`}
\def\@stressit{\futurelet\chartest\@stresschar }

\def\@stresschar#1{%
  \ifx #1y\def\result{\futurelet\chartest\@yligature}%
  \else \ifx #1Y\def\result{\futurelet\chartest\@Yligature}%
  \else \ifx\chartest\@testgrave \def\result{\accent"26 }%
  \else \def\result{\accent"26 #1}%
  \fi \fi \fi
  \result }

\def\@yligature{%
  \ifx a\chartest \def\result{\accent"26 \char"1F \@gobble}%
  \else \ifx u\chartest \def\result{\accent"26 \char"18 \@gobble}%
  \else \def\result{\accent"26 y}%
  \fi \fi
  \result }

\def\@Yligature{%
  \ifx a\chartest \def\result{\accent"26 \char"17 \@gobble}%
  \else \ifx A\chartest \def\result{\accent"26 \char"17 \@gobble}%
  \else \ifx u\chartest \def\result{\accent"26 \char"10 \@gobble}%
  \else \ifx U\chartest \def\result{\accent"26 \char"10 \@gobble}%
  \else \def\result{\accent"26 Y}%
  \fi \fi \fi \fi
  \result }

\def\!{\ifmmode \mskip-\thinmuskip \fi}


\def\cyracc{%
  \def\cydot{{\kern0pt}}%
  \def\cprime{\char"7E }\def\Cprime{\char"5E }%
  \def\cdprime{\char"7F }\def\Cdprime{\char"5F }%
  \def\dbar{dj}\def\Dbar{Dj}%
  \def\dz{\char"1E }\def\Dz{\char"16 }%
  \def\dzh{\char"0A }\def\Dzh{\char"02 }%
  \def\'##1{\if c##1\char"0F %
    \else \if C##1\char"07 %
    \else \accent"26 ##1\fi \fi }%
  \def\`##1{\if e##1\char"0B %
    \else \if E##1\char"03 %
    \else \errmessage{accent \string\` not defined in cyrillic}%
        ##1\fi \fi }%
  \def\=##1{\if e##1\char"0D %
    \else \if E##1\char"05 %
    \else \if \i##1\char"0C %
    \else \if I##1\char"04 %
    \else \errmessage{accent \string\= not defined in cyrillic}%
        ##1\fi \fi \fi \fi }%
  \def\u##1{\if \i##1\accent"24 i%
    \else \accent"24 ##1\fi }%
  \def\"##1{\if \i##1\accent"20 \char"3D %
    \else \if I##1\accent"20 \char"04 %
    \else \accent"20 ##1\fi \fi }%
  \def\!{\ifmmode \def\result{\mskip-\thinmuskip}%
    \else \def\result{\@stressit}\fi \result}}


\def\keep@cyracc{\let\cyr=\relax \let\i=\relax
        \let\ubar=\relax \let\cydot=\relax
        \let\cprime=\relax \let\Cprime=\relax
        \let\cdprime=\relax \let\Cdprime=\relax
        \let\dbar=\relax \let\Dbar=\relax
        \let\dz=\relax \let\Dz=\relax
        \let\dzh=\relax \let\Dzh=\relax
        \let\'=\relax \let\`=\relax \let\==\relax
        \let\u=\relax \let\"=\relax \let\!=\relax }

\catcode`\@=\tempcat
\newfam\cyrfam
\font\twelvecyr=wncyr10 scaled 1200
\def\cyr{\fam\cyrfam\twelvecyr\cyracc}

\newcount\EQNO      \EQNO=0
\newcount\FIGNO     \FIGNO=0
\newcount\REFNO     \REFNO=0
\newcount\SECNO     \SECNO=0
\newcount\SUBSECNO  \SUBSECNO=0
\newcount\FOOTNO    \FOOTNO=0
\newbox\FIGBOX      \setbox\FIGBOX=\vbox{}
\newbox\REFBOX      \newbox\REFBOXTMP
\setbox\REFBOX=\vbox{\bigskip\centerline{\bf REFERENCES}\smallskip}
\newbox\Partialpage \newdimen\Mark
\newdimen\REFSIZE   \REFSIZE=\vsize

\expandafter\ifx\csname normal\endcsname\relax\def\normal{\null}\fi

\def\MultiRef#1{\global\advance\REFNO by 1 \nobreak\the\REFNO%
    \global\setbox\REFBOX=\vbox{\unvcopy\REFBOX\normal
      \smallskip\item{\the\REFNO .~}#1}%
    \gdef\label##1{\xdef##1{\nobreak[\the\REFNO]}}%
    \gdef\Label##1{\xdef##1{\nobreak\the\REFNO}}}
\def\NoRef#1{\global\advance\REFNO by 1%
    \global\setbox\REFBOX=\vbox{\unvcopy\REFBOX\normal
      \smallskip\item{\the\REFNO .~}#1}%
    \gdef\label##1{\xdef##1{\nobreak[\the\REFNO]}}}
\def\Eqno{\global\advance\EQNO by 1 \eqno(\the\EQNO)%
    \gdef\label##1{\xdef##1{\nobreak(\the\EQNO)}}}
\def\Eqalignno{\global\advance\EQNO by 1 &(\the\EQNO)%
    \gdef\label##1{\xdef##1{\nobreak(\the\EQNO)}}}
\def\Fig#1{\global\advance\FIGNO by 1 Figure~\the\FIGNO%
    \global\setbox\FIGBOX=\vbox{\unvcopy\FIGBOX
      \narrower\smallskip\item{\bf Figure \the\FIGNO~~}#1}}
\def\Ref#1{\global\advance\REFNO by 1 \nobreak[\the\REFNO]%
    \global\setbox\REFBOX=\vbox{\unvcopy\REFBOX\normal
      \smallskip\item{\the\REFNO .~}#1}%
    \gdef\label##1{\xdef##1{\nobreak[\the\REFNO]}}}
\def\Section#1{\SUBSECNO=0\advance\SECNO by 1
    \bigskip\leftline{\bf \the\SECNO .\ #1}\nobreak}
\def\Subsection#1{\advance\SUBSECNO by 1
    \medskip\leftline{\bf \ifcase\SUBSECNO\or
    a\or b\or c\or d\or e\or f\or g\or h\or i\or j\or k\or l\or m\or n\fi
    )\ #1}\nobreak}
\def\Footnote#1{\global\advance\FOOTNO by 1 
    \footnote{\nobreak$\>\!{}^{\the\FOOTNO}\>\!$}{#1}
    \gdef\label##1{\xdef##1{$\>\!{}^{\the\FOOTNO}\>\!$}}}
\def\SameFootnote{$\>\!{}^{\the\FOOTNO}\>\!$}

\def\References{
     {\output={\global\setbox\Partialpage=\vbox{\unvbox255}}\eject}
     \Mark=\vsize
     \ifdim \ht\Partialpage > 0pt
          \advance\Mark by -\ht\Partialpage
          \advance\Mark by -0.1in
          \vbox{\unvbox\Partialpage}
     \fi
     \setbox\REFBOX=\vbox{\unvcopy\REFBOX\vfill\eject}
     \ifdim \ht\REFBOX > \Mark
          \setbox\REFBOXTMP=\vsplit\REFBOX to \Mark
          \vfill\copy\REFBOXTMP
     \else
          \REFSIZE=\Mark
     \fi
     \loop
          \setbox\REFBOXTMP=\vsplit\REFBOX to \REFSIZE
          \copy\REFBOXTMP
          \ifdim \ht\REFBOX > 0pt
     \repeat}
\def\NewRefPage{\setbox\REFBOX=\vbox{\unvcopy\REFBOX\vfill\eject}}


\newcount\LEMMA	  \LEMMA=0
\newcount\THM     \THM=0
\def\Lemma#1{\global\advance\LEMMA by 1\smallskip
    {\narrower\narrower\narrower\item{\bf Lemma~\the\LEMMA:}
    {\it #1}\smallskip}}
\def\Theorem#1{\global\advance\THM by 1\smallskip
    {\narrower\narrower\narrower\item{\bf Theorem~\the\THM:}
    {\it #1}\smallskip}}

\input epsf
\def\Fig#1#2#3#4#5{\global\advance\FIGNO by 1 Figure~\the\FIGNO#5
    \topinsert
    \centerline{\epsfysize=#4\epsffile[#3]{#2}}
    {\bigskip\hsize=5.5in\hskip\parindent
     \vbox{\small\item{\bf Figure~\the\FIGNO:}{#1}}}
    \bigskip\endinsert}

\def\RR{{\Bbb R}}
\def\CC{{\Bbb C}}

\def\OO{{\Bbb O}}
\def\AA{{\cal A}}
\def\BB{{\cal B}}

\def\bar{\overline}

\def\Tr{{\rm tr\,}}
\def\Re{{\rm Re}}
\def\Im{{\rm Im}}


\footline={
	\ifnum\pageno>1
		{\hss\rm \folio \hss}
	\else
		{\small\hss Copyright \copyright\ 2000 by Tevian Dray\hss}
	\fi}

\rightline{
	\small{\tt http://xxx.lanl.gov/abs/math/0010255}\hfill 11 October 2000}
\bigskip

\centerline{\large\bf SOME PROPERTIES OF {\largebold$3\times3$}}
\smallskip
\centerline{\large\bf OCTONIONIC HERMITIAN MATRICES}
\smallskip
\centerline{\large\bf WITH NON-REAL EIGENVALUES}
\bigskip

\centerline{Tevian Dray}
\centerline{\it Department of Mathematics, Oregon State University,
		Corvallis, OR  97331, USA}
\centerline{\tt tevian{\rm @}math.orst.edu}
\medskip
\centerline{Jason Janesky}
\centerline{\it Department of Physics, Oregon State University,
		Corvallis, OR  97331, USA
\Footnote{Present address:
	Phoenix Corporate Research Laboratories,
	Motorola Inc.,
	Tempe, AZ  85284, USA,
	{\tt r47569{\rm @}email.sps.mot.com}
}}

\medskip
\centerline{Corinne A. Manogue}
\centerline{\it Department of Physics, Oregon State University,
		Corvallis, OR  97331, USA}
\centerline{\tt corinne{\rm @}physics.orst.edu}

\bigskip\bigskip
\centerline{\bf ABSTRACT}
\midinsert
\narrower\narrower\noindent
We discuss our preliminary attempts to extend previous work on $2\times2$
Hermitian octonionic matrices with non-real eigenvalues to the $3\times3$
case.
\endinsert
\bigskip

\Section{INTRODUCTION}

In previous work
[\MultiRef{Tevian Dray and Corinne A. Manogue,
{\it The Octonionic Eigenvalue Problem},
Adv.\ Appl.\ Clifford Algebras {\bf 8}, 341--364 (1998).}%
\label\Eigen\Label\EIGEN
,\MultiRef{Tevian Dray and Corinne A. Manogue,
{\it Finding Octonionic Eigenvectors Using {\sl Mathematica}},
Comput.\ Phys.\ Comm.\ {\bf 115}, 536--547 (1998).}\label\Find\Label\FIND
; see also
\MultiRef{Susumu Okubo,
{\it Eigenvalue Problem for Symmetric $3\times3$ Octonionic Matrix},
Adv.\ Appl.\ Clifford Algebras {\bf 9}, 131--176 (1999).}]\label\Okubo
, we considered the real eigenvalue problem for $2\times2$ and $3\times3$
Hermitian matrices over the octonions $\OO$.  The $2\times2$ case corresponds
closely to the standard, complex eigenvalue problem, since any $2\times2$
octonionic Hermitian matrix lies in a complex subalgebra $\CC\subset\OO$.  The
$3\times3$ case requires considerable care, resulting in some changes in the
expected results.  However, we also showed in \Eigen\ that there are
octonionic Hermitian matrices which admit eigenvalues which are not real, and
a complete treatment of the $2\times2$ case was given in
\Ref{Tevian Dray, Jason Janesky, and Corinne A. Manogue,
{\it Octonionic Hermitian Matrices with Non-Real Eigenvalues},
Adv.\ Appl.\ Clifford Algebras (submitted).
\hfill{\small({\tt math/0006069})}}\label\Nonreal
.  Here, we discuss our preliminary results for the $3\times3$ case.

Although we are able to obtain a 3rd-order characteristic equation for the
(right) eigenvalues in the $3\times3$ case, we have not been able to solve
this equation, nor have we been able to extend our orthonormality results
[\EIGEN,\FIND] from the real case.  We therefore discuss several illustrative
examples and make some conjectures regarding more general results.

The $3\times3$ case is of particular interest mathematically because it
corresponds to the exceptional Jordan algebra, also known as the Albert
algebra.  There have been numerous attempts to use this algebra to describe
quantum physics, which was in fact Jordan's original motivation.  More
recently, Schray 
[\MultiRef{J\"org Schray,
{\bf Octonions and Supersymmetry},
Ph.D.\ thesis, Department of Physics, Oregon State University, 1994.},%
\MultiRef{J\"org Schray,
{\it The General Classical Solution of the Superparticle},
Class.\ Quant.\ Grav.\ {\bf 13}, 27 (1996)}] \label\Schray\Label\SCHRAY
has shown how to use the exceptional Jordan algebra to give
an elegant description of the superparticle, which we have been attempting to
extend to the superstring.  Our dimensional reduction scheme extends naturally
to this case
\Ref{Tevian Dray and Corinne A. Manogue,
{\it The Exceptional Jordan Eigenvalue Problem},
Internat.\ J.\ Theoret.\ Phys.\ {\bf 38}, 2901--2916 (1999).}\label\Jordan
, and we believe it is the natural language to describe
the fundamental particles of nature.

The paper is organized as follows.  In Section 2 we briefly review the
properties of octonions.  We then consider $3\times3$ octonionic Hermitian
matrices, deriving a characteristic equation for the eigenvalues in Section 3,
and considering several examples in Section 4.  Finally, in Section 5 we
discuss our results.  Some parts of this presentation have appeared in our
previous work.

\Section{OCTONIONS}

We summarize here only the essential properties of the octonions $\OO$.  For a
more detailed introduction, see \Eigen\ or
[\MultiRef{Feza G\"ursey and Chia-Hsiung Tze,
{\bf On the Role of Division, Jordan, and Related Algebras in Particle Physics},
World Scientific, Singapore, 1996.}%
,\MultiRef{S. Okubo,
{\bf Introduction to Octonion and Other Non-Associative Algebras in Physics},
Cambridge University Press, Cambridge, 1995.}]%
.%

The octonions $\OO$ are the nonassociative, noncommutative, normed division
algebra over the reals.  In terms of a natural basis, an octonion $a$ can be
written
$$a = \sum\limits_{q=1}^8 a^q e_q \Eqno$$
where the coefficients $a^q$ are real, and where the basis vectors satisfy
$e_1=1$ and
$$e_q^2 = -1 \qquad (q=2,...,8) \Eqno$$
The multiplication table is
conveniently encoded in the 7-point projective plane, shown in
\Fig{The representation of the octonionic multiplication table using the
7-point projective plane, where we have used the conventional names
$\{i,j,k,k\ell,j\ell,i\ell,\ell\}$ for $\{e_2,...,e_8\}$.  Each of the 7
oriented lines gives a quaternionic triple.}
{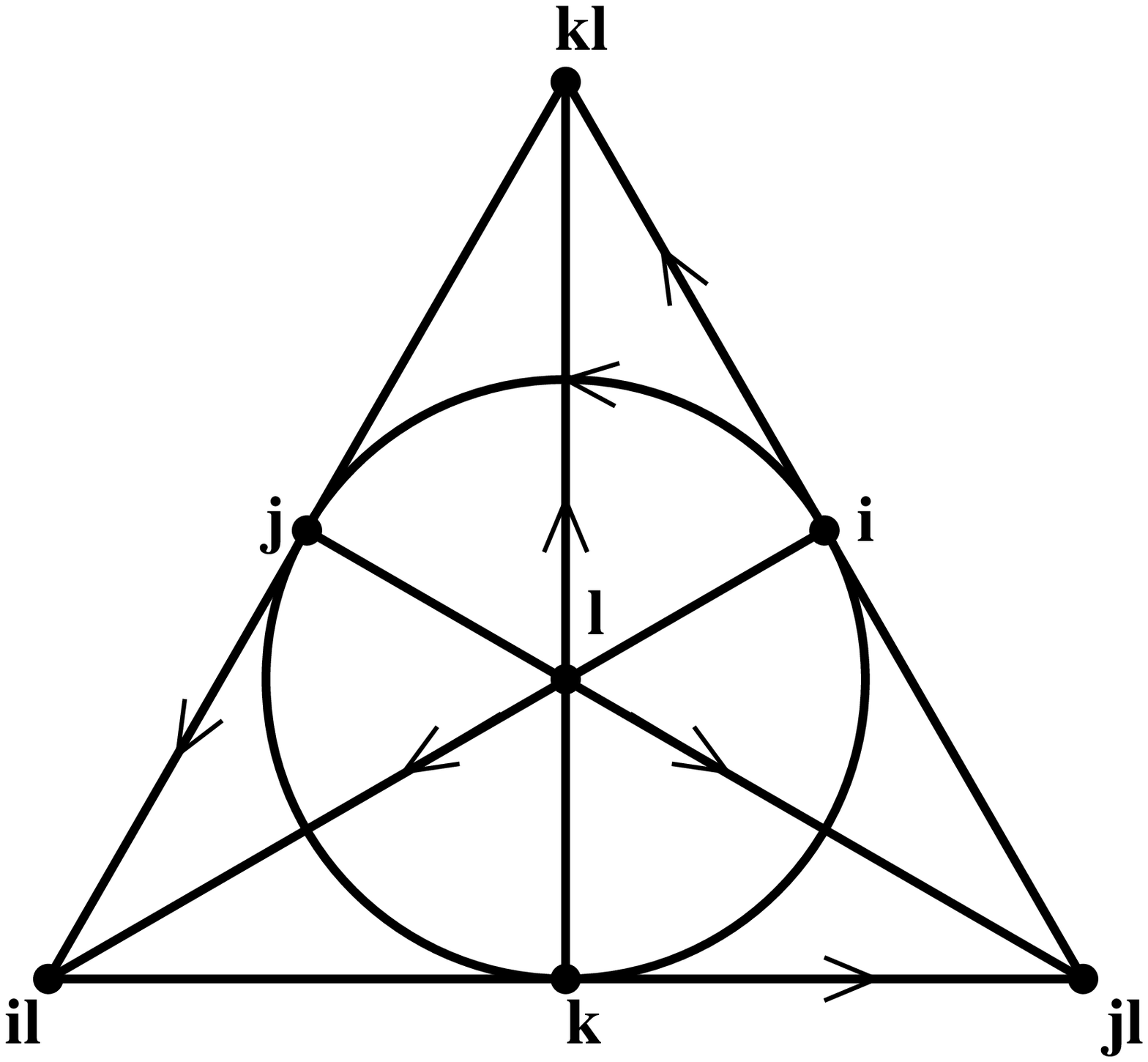}{68 168 543 614}{3in}
{.  The product of any two imaginary units is given by the third unit on the
unique line connecting them, with the sign determined by the relative
orientation.}

{\it Octonionic conjugation\/} is given by reversing the sign of the
imaginary basis units
$$\bar a = a^1 e_1 - \sum\limits_{q=2}^8 a^q e_q \Eqno$$
Conjugation is an antiautomorphism, since it satisfies
$$\bar{ab} = \bar{b} \> \bar{a}$$
The real and imaginary parts of an octonion $a$ are given by
$$\Re(a) = {1\over2} (a+\bar{a})
  \qquad\qquad
  \Im(a) = {1\over2} (a-\bar{a})
  \Eqno$$

The {\it inner product\/} on $\OO$ is the one inherited from $\RR^8$, namely
$$a \cdot b = \sum_q a^q b^q \Eqno$$
which can be rewritten as
$$a \cdot b
  = {1\over2} (a \bar{b} + b \bar{a})
  = {1\over2} (\bar{b} a + \bar{a} b)
  \Eqno$$
and which satisfies the identities
$$\eqalignno{
a \cdot (xb) &= b \cdot (\bar{x}a) \Eqalignno\label\DotID \cr
(ax) \cdot (bx) &= |x|^2 \> a \cdot b \Eqalignno\label\DotIDii
  }$$
for any $a,b,x\in\OO$.  The {\it norm\/} of an octonion is just
$$|a| = \sqrt{a \bar{a}} = \sqrt{a \cdot a} \Eqno$$
which satisfies the defining property of a normed division algebra, namely
$$|ab| = |a| |b| \Eqno$$

The {\it associator\/} of three octonions is
$$[a,b,c] = (ab)c - a(bc) \Eqno$$
which is totally antisymmetric in its arguments, has no real part, and changes
sign if any one of its arguments is replaced by its octonionic conjugate.
Although the associator does not vanish in general, the octonions do satisfy a
weak form of associativity known as {\it alternativity}, namely
$$[b,a,a]=0=[b,a,\bar{a}] \Eqno$$
The underlying reason for alternativity is Artin's Theorem
[\MultiRef{
Richard D. Schafer,
{\bf An Introduction to Nonassociative Algebras},
Academic Press, New York, 1966; reprinted by Dover, Mineola NY, 1995.},%
\MultiRef{Emil Artin,
{\bf Geometric Algebra},
Interscience Publishers, New York, 1957; reprinted by
John Wiley \& Sons, New York, 1988.}\label\Artin
]%
, which states that that any two octonions lie in a quaternionic subalgebra of
$\OO$, so that any product containing only two octonionic directions is
associative.  We will also have use for the associator identity
$$[a,b,c]d + a[b,c,d] = [ab,c,d] - [a,bc,d] + [a,b,cd] \Eqno$$\label\Assoc
for any $a,b,c,d\in\OO$, which is proved by writing out all the terms.

\Section{\Bold $3\times3$ OCTONIONIC HERMITIAN MATRICES}

In this section, we derive a characteristic equation for the (right)
eigenvalues of $\AA$, which reduces to that of \Eigen\ for real eigenvalues.
Unfortunately, we have been unable to solve this equation when the eigenvalues
are not real, so that we have also been unable to investigate orthogonality
and decomposition results analogous to those for real eigenvalues.  We discuss
this further in the next section, where we study several examples with
intriguing properties.

\Subsection{Jordan matrices}

The $3\times3$ octonionic Hermitian matrices, henceforth referred to as
{\it Jordan matrices}, form the exceptional Jordan algebra (also called the
Albert algebra) under the Jordan product
\Footnote{The $2\times2$ octonionic Hermitian matrices form a special Jordan
algebra since they are alternative
\Ref{Nathan Jacobson,
{\bf Structure and Representations of Jordan Algebras},
Amer.\  Math.\  Soc.\ Colloq.\ Publ.\ {\bf 39}, 
American Mathematical Society, Providence, 1968.}\label\Jacobson
.}
$$\AA \circ \BB := {1\over2} (\AA\BB + \BB\AA) \Eqno$$
which is commutative, but not associative.  A special case of this is
$$\AA^2 \equiv \AA \circ \AA \Eqno$$
and we {\it define}
$$\AA^3 := \AA^2 \circ \AA = \AA \circ \AA^2 \Eqno$$

Remarkably, with these definitions, Jordan matrices satisfy the usual
characteristic equation
\Ref{Hans Freudenthal,
{\it Zur Ebenen Oktavengeometrie},
Proc.\ Kon.\ Ned.\ Akad.\ Wet.\ {\bf A56}, 195--200 (1953).}
$$\AA^3 - (\Tr \AA) \, \AA^2 + \sigma(\AA) \, \AA - (\det \AA) \, I
  = 0 \Eqno$$\label\Char
where $\sigma(\AA)$ is defined by
$$\sigma(\AA) := {1\over2} \left( (\Tr \AA)^2 - \Tr (\AA^2) \right) \Eqno$$
and where the determinant of $\AA$ is defined abstractly in terms of the
Freudenthal product.~%
\Footnote{
The Freudenthal product of two Jordan matrices $\AA$ and $\BB$ is given by
\Ref{P. Jordan, J. von Neumann, and E. Wigner,
Ann.\ Math.\ {\bf 35}, 29 (1934);
\hfill\break
H. Freudenthal, Adv.\ Math.\ {\bf 1}, 145 (1964).}
$$\AA*\BB = \AA \circ \BB - {1\over2} \Big(\AA\,\Tr(\BB)+\BB\,\Tr(\AA)\Big)
		+ {1\over2} \Big(\Tr(\AA)\,\Tr(\BB)-\Tr(\AA\circ \BB)\Big)$$
The determinant can then be defined as
$$\det(\AA) = {1\over3} \, \Tr \Big( (\AA*\AA) \circ \AA \Big)$$
}
Concretely, if
$$\AA = \pmatrix{p& a& \bar{b}\cr \bar{a}& m& c\cr b& \bar{c}& n\cr} \Eqno$$
\label\Three
with $p,m,n\in\RR$ and $a,b,c\in\OO$ then
$$\eqalign{
\Tr \AA &= p + m + n \cr
\sigma(\AA) &= pm + pn + mn - |a|^2 - |b|^2 - |c|^2 \cr
\det \AA &= pmn + b(ac) + \bar{b(ac)} - n|a|^2 - m|b|^2 - p|c|^2 \cr
  }\Eqno$$\label\ThreeEq

As shown originally by Ogievetsky
\Ref{
{\cyr O. V. Ogievetski\u\i},
{\cyr Kharakteristicheskoe U$\!$ravnenie dlya Matrits ${3\times3}$
		nad Oktavami}, 
Uspekhi Mat.\ Nauk {\bf 36}, 197--198 (1981); translated in:
O. V. Ogievetskii,
{\it The Characteristic Equation for ${3\times3}$ Matrices
	over Octaves},
Russian Math.\ Surveys {\bf 36}, 189--190 (1981).}%
, $\AA$, has 6, rather than 3, real eigenvalues, which furthermore fail to
satisfy the characteristic equation \Char.  A complete, computer-assisted
\Find\ treatment of this case was given in \Eigen, and a somewhat more general
analytic treatment was later given by Okubo \Okubo.  As shown there, the
eigenvalues naturally belong to 2 distinct families, each containing 3 real
eigenvalues. Furthermore, within each family, the corresponding eigenvectors
lead to a decomposition of the form
$$\AA
  = \sum_{\alpha=1}^3 \lambda_\alpha \left( v_\alpha v_\alpha^\dagger \right)
  \Eqno$$\label\DecompGen
Furthermore, eigenvectors $v_\alpha$ corresponding to different eigenvalues
are automatically orthogonal in the generalized sense
$$(v v^\dagger) \, w = 0 \Eqno$$\label\Ortho

\Subsection{Characteristic Equation}

Set
$$v=\pmatrix{x\cr y\cr z\cr} \Eqno$$\label\vIII
Then the (right) eigenvalue problem
$$A v = v\lambda \Eqno$$\label\Master
becomes
$$\eqalignno{
  x \, (\lambda-p) &= ay + \bar{b}z \Eqalignno\cr \label\EI
  y \, (\lambda-m) &= cz + \bar{a}x \Eqalignno\cr \label\EII
  z \, (\lambda-n) &= bx + \bar{c}y \Eqalignno \label\EIII
  }$$
As shown in \Eigen, \Master\ admits solutions for which $\lambda$ is not real;
further examples are given in Section 4.
\Footnote{The quite different {\it Jordan\/} eigenvalue problem in the
$3\times3$ case admits only real eigenvalues and was discussed in \Jordan.}
Multiplying \EI\ on the right by $(\lambda-m)$ leads to
$$(ay) (\lambda-m) = x (\lambda-p) (\lambda-m) - (\bar{b}z) (\lambda-m) \Eqno$$
whereas multiplying \EII\ on the left by $\bar{a}$
$$a \Big(y(\lambda-m)\Big) = a (cz) + |a|^2 x \Eqno$$
Subtracting these 2 equations immediately yields
$$[a,y,\lambda]
  = x \left( (\lambda-p)(\lambda-m) - |a|^2 \right)
    - (\bar{b}z) (\lambda-m) - a (cz)
\Eqno$$\label\FI
Similarly, multiplying \EII\ on the right by $(\lambda-p)$ and \EI\ on the
left by $\bar{a}$ (or using symmetry) leads to
$$[\bar{a},x,\lambda]
  = y \left( (\lambda-p)(\lambda-m) - |a|^2 \right)
    - (cz) (\lambda-p) - \bar{a} (\bar{b}z)
\Eqno$$\label\FII

We plan to multiply \FI\ by $b$ on the left, \FII\ by $\bar{c}$ on the left,
add, and then use \EIII.  Before doing so, we first use \Assoc\ to write
$$\eqalign{
  [b,x,(\lambda-p)(\lambda-m)]
  &= [b,x(\lambda-p),(\lambda-m)] - [bx,(\lambda-p),(\lambda-m)] \cr
  &\qquad + b \, [x,(\lambda-p),(\lambda-m)] + [b,x,(\lambda-p)] \, (\lambda-m)
	\cr
  &= [b,x(\lambda-p),\lambda] + [b,x,\lambda] \, (\lambda-m) \cr
  &= [b,ay+\bar{b}z,\lambda] + [b,x,\lambda] \, (\lambda-m) \cr
  }\Eqno$$
as well as
$$\eqalign{
  [b,\bar{b}z,\lambda]
  &= [b\bar{b},z,\lambda] + [b,\bar{b},z\lambda]
     - b \, [\bar{b},z,\lambda] - [b,\bar{b},z] \, \lambda \cr
  &= - b \, [\bar{b},z,\lambda] \cr
  }\Eqno$$
Thus, returning to \FI, we obtain
$$\eqalign{
  b \, [a,y,\lambda]
  &= (bx) \left( (\lambda-p)(\lambda-m) - |a|^2 \right)
     - |b|^2 z (\lambda-m) - b \Bigl( a (cz) \Bigr) \cr
  &\qquad - [b,x,(\lambda-p)(\lambda-m)] - b \, [\bar{b},z,\lambda] \cr
  &= (bx) \left( (\lambda-p)(\lambda-m) - |a|^2 \right)
     - |b|^2 z (\lambda-m) - b \Bigl( a (cz) \Bigr) \cr
  &\qquad - [b,ay,\lambda] - [b,x,\lambda] \, (\lambda-m) \cr
  }\Eqno$$
Similarly, \FII\ becomes
$$\eqalign{
  \bar{c} \, [\bar{a},x,\lambda]
  &= (\bar{c}y) \left( (\lambda-p)(\lambda-m) - |a|^2 \right)
     - |c|^2 z (\lambda-p) - \bar{c} \Bigl( \bar{a} (\bar{b}z) \Bigr) \cr
  &\qquad - [\bar{c},\bar{a}x,\lambda] - [\bar{c},y,\lambda] \, (\lambda-p) \cr
  }\Eqno$$
Adding these last 2 equations, we obtain
$$\eqalign{
  b \, [a,y,\lambda] + \bar{c} \, [\bar{a},x,\lambda]
  &= (bx+\bar{c}y) \left( (\lambda-p)(\lambda-m) - |a|^2 \right)
     - |b|^2 z (\lambda-m) - |c|^2 z (\lambda-p) \cr
  & \qquad
     - b \Bigl( a (cz) \Bigr) - \bar{c} \Bigl( \bar{a} (\bar{b}z) \Bigr) \cr
  &\qquad - [\bar{c},\bar{a}x,\lambda] - [\bar{c},y,\lambda] \, (\lambda-p)
          - [b,ay,\lambda] - [b,x,\lambda] \, (\lambda-m) \cr
  }\Eqno$$
Finally, using \EIII, factoring out $z$, and rearranging terms leads to the
generalized characteristic equation in the form
\Footnote{The extra terms multiplying $z$ on the right-hand-side come from
$\det \AA$.}
$$\eqalign{
  z \left(
	\lambda^3 - (\Tr \AA) \, \lambda^2 + \sigma(\AA) \, \lambda - \det \AA
	\right)
  &= b \Bigl( a (cz) \Bigr) + \bar{c} \Bigl( \bar{a} (\bar{b}z)	\Bigr) 
     - \Bigl( b(ac) + (\bar{c}\,\bar{a})\bar{b} \Bigr) z \cr
  &\qquad + b \, [a,y,\lambda] + [b,ay,\lambda]
          + [b,x,\lambda] \, (\lambda-m) \cr 
  &\qquad + \bar{c} \, [\bar{a},x,\lambda] + [\bar{c},\bar{a}x,\lambda]
          + [\bar{c},y,\lambda] \, (\lambda-p) \cr
  }\Eqno$$\label\CharIII

If $\lambda$ is real, all the associators on the right-hand-side vanish, and
we recover the generalized characteristic equation given in \Eigen.  The
requirement in that case that the right-hand-side be a real multiple of $z$
(since the left-hand-side is) then constrains $z$, resulting in precisely 2
values for that real multiple, and reducing \CharIII\ to 2 cubic equations,
one for each family of real eigenvalues.

While we find the form of \CharIII\ attractive, as there are no extraneous
terms involving both $z$ and $\lambda$, we have so far been unable to further
simplify \CharIII\ when $\lambda$ is not real.

\Subsection{Alternate Approach}

We briefly describe another possible approach to finding the eigenvalues,
which relies on the associator identity
$$[v^\dagger,v,\lambda]
  := (v^\dagger v) \lambda - v^\dagger (v\lambda) \equiv 0 \Eqno$$\label\APP
which follows for {\it any} octonionic vector $v$ and $\lambda\in\OO$ by
alternativity, and which is further discussed in the Appendix of \Nonreal.  If
$v$ is a normalized right eigenvector of $\AA$ with eigenvalue $\lambda$, then
$$v^\dagger (\AA v) = v^\dagger (v\lambda) = (v^\dagger v) \lambda = \lambda
  \Eqno$$\label\Jason
which yields an equation for $\lambda$ in terms of $\AA$ and the components of
$v$.  A similar construction using the associator
$$[v^\dagger,\AA,v]
  := (v^\dagger \AA) v - v^\dagger (\AA v)
  = (\AA v)^\dagger v - v^\dagger (\AA v)
  = \left( v^\dagger (\AA v) \right)^\dagger - v^\dagger (\AA v)
   \Eqno$$
leads for normalized eigenvectors to
$$[v^\dagger,\AA,v]
  = \left( v^\dagger (v\lambda) \right)^\dagger - v^\dagger (v\lambda)
  = \left( (v^\dagger v) \lambda) \right)^\dagger - (v^\dagger v) \lambda
  = \bar\lambda-\lambda = -2 \, \Im(\lambda)
\Eqno$$

Inserting \Three\ and \vIII\ into \Jason\ leads, after minor rearrangement
using associators and \EIII, to
$$\lambda
  = {p|x|^2 + m|y|^2 - n|z|^2 + 2 x \cdot (ay) \over |x|^2 + |y|^2 - |z|^2}
    + {[x,a,y] + [z,b,x] + [y,c,z] \over |x|^2 + |y|^2 + |z|^2}
  \Eqno$$\label\LambdaEq
which gives explicit expressions for the real and imaginary parts of
$\lambda$.  The first term can be rewritten using cyclic permutations of
$\{x,y,z\}$ (and $\{a,c,b\}$), and the resulting expressions set equal to
obtain
$$\Re(\lambda) = {x \cdot (ay) + z \cdot (bx) + p |x|^2 \over |x|^2} \Eqno$$
and similar expressions obtained by cyclic permutation.
\Footnote{These expressions can also be obtained directly from \EI--\EIII\ by
taking the dot product with $x$, $y$, $z$, respectively, and using \DotID.}
Finally, if $v$ is normalized, the imaginary part of \LambdaEq\ reduces to
$$\Im(\lambda) = [x,a,y] + [z,b,x] + [y,c,z] \Eqno$$
We had hoped to use these various expressions to impose conditions on $\AA$
which would in turn enable us to solve for $\lambda$, but have not yet found a
way to do so.

\Section{EXAMPLES}

Without being able to solve (some version of) the characteristic equation in
the $3\times3$ case, it is not possible in general to determine all the
(non-real) eigenvalues of a given Hermitian octonionic matrix.  It is
therefore instructive to consider several explicit examples.

\Subsection{Example 1}

Consider the matrix
$$\BB = \pmatrix{~~p & ~~iq & kqs \cr \noalign{\smallskip}
		-iq & ~~p & jq\cr \noalign{\smallskip}
		-kqs & -jq & p \cr}
 \Eqno$$
where
$$s=\cos\theta+k\ell\sin\theta \Eqno$$\label\Seq
Note that $\BB$ is quaternionic if $\theta=0$.

The real eigenvalues of $\BB$, and corresponding orthonormal bases of
eigenvectors, were given in \Find.  But $\BB$ also admits eigenvectors with
eigenvalues which are not real.  For instance:
$$\eqalign{
\lambda_{\hat u} = p \pm q\bar{s}: \quad&
	\hat u_\pm = \pmatrix{i\cr 0\cr j\cr} S_\pm \cr
\noalign{\smallskip}
\lambda_{\hat v} = p \pm q\bar{s}: \quad&
	\hat v_\pm = \pmatrix{j\cr 2ks\cr i\cr} S_\pm \cr
\noalign{\smallskip}
\lambda_{\hat w} = p \mp 2q\bar{s}: \quad&
	\hat w_\pm = \pmatrix{~j\cr -ks\cr ~i\cr} S_\pm \cr
  }\Eqno$$
where
$$S_\pm = \cases{-k\ell\cr ~~1\cr} \Eqno$$
These eigenvectors and eigenvalues reduce to the ones given in \Find\ when
$\theta\to0$.  Somewhat surprisingly, these eigenvectors (when normalized)
yield a decomposition of the form~\DecompGen.  Remarkably, they also yield a
decomposition of the form
$$\AA
  = \sum_{\alpha=1}^3 \left( v_\alpha \lambda_\alpha \right) v_\alpha^\dagger
  \Eqno$$\label\DecompSp

\goodbreak
We now describe some further properties of these eigenvectors.  Each
eigenspace with eigenvalues $\lambda_{\hat w} = p \mp 2q\bar{s}$ is
1-dimensional (over \RR), so that the eigenvectors $\hat{w}_\pm$ are
essentially unique.  By contrast, the eigenspaces with eigenvalues
$\lambda_{\hat u} =
\lambda_{\hat v} = p \pm q\bar{s}$ are 5-dimensional.  Interestingly, though,
for any given eigenvector such as $\hat{u}_\pm$, there is again an essentially
unique eigenvector, in this case $\hat{v}_\pm$, which is orthogonal to it.
Here ``essentially unique'' means unique up to a real multiplicative factor,
and orthogonality can be defined either as $v^\dagger w = 0$ or as
$(vv^\dagger) w = 0$; these turn out to be equivalent in this case.

There are also additional eigenvectors with eigenvalues of the form
$$\lambda = (p+\rho) - \beta \, k\ell \Eqno$$\label\Lform
where $\beta,\rho\in\RR$.  These must satisfy
$$32\beta^2 = \left(\sqrt{32\rho^2-7q^2}-5q\right)
		\left(11q-\sqrt{32\rho^2-7q^2}\right) \Eqno$$
and therefore only exist provided that
$$q^2 \le \rho^2 \le 4 q^2 \Eqno$$
In the case of equality, we recover the real eigenvalues given in \Find.  For
each admissible $\rho$ and each of the two corresponding choices for
$\beta$, the eigenspace is 3-dimensional.  We have not explored the properties
of these eigenvectors in depth.

All the eigenvalues discussed above have the form \Lform.  While we suspect
that there are no others, we have not been able to prove that this is the
case.

\Subsection{Example 2}

A related example is given by the matrix
$$\hat{\BB} = \pmatrix{~~p & ~~qi & {q\over2}ks \cr \noalign{\smallskip}
		-qi & ~~p & {q\over2}j\cr \noalign{\smallskip}
		-{q\over2}ks & -{q\over2}j & p \cr}
 \Eqno$$
with $s$ again given by \Seq.  We choose $\theta$ such that
$$s= {\sqrt{5}\over3} - {2\over3} k\ell \Eqno$$
resulting in
$$\hat{\BB} = \pmatrix{~~p & ~~qi & {q\over6}(\sqrt{5}k+2l) \cr
		\noalign{\smallskip}
		-qi & ~~p & {q\over2}j\cr
		\noalign{\smallskip}
		-{q\over6}(\sqrt{5}k+2l) & -{q\over2}j & p \cr}
 \Eqno$$
The 2 families of real eigenvalues of $\hat{\BB}$ turn out to be
$\{p\pm q, p\mp{q\over2}(1+\sqrt{3}), p\mp{q\over2}(1-{\sqrt{3}\over2})\}$.
Some eigenvectors for $\hat{\BB}$ corresponding to eigenvalues which are not
real are
$$\eqalign{
\lambda_{u_1} = (p+{\sqrt{5}\over2}q) - {q\over2}k\ell: \quad&
	u_1 = \pmatrix{3k\cr \sqrt{5}\,j-2\,i\ell\cr 1+\sqrt{5}\,k\ell\cr} \cr
\noalign{\smallskip}
\lambda_{u_2} = (p+{\sqrt{5}\over2}q) + {q\over2}k\ell: \quad&
	u_2 = \pmatrix{\sqrt{5}\,k+2\ell\cr 3j\cr \sqrt{5}-k\ell\cr} \cr
\noalign{\smallskip}
\lambda_{v_1} = (p-{\sqrt{5}\over3}q) + {2q\over3}k\ell: \quad&
	v_1 = \pmatrix{\sqrt{5}\,j-2\,i\ell\cr 3k\cr 0\cr} \cr
\noalign{\smallskip}
\lambda_{v_2} = (p-{\sqrt{5}\over3}q) - {2q\over3}k\ell: \quad&
	v_2 = \pmatrix{3j\cr \sqrt{5}\,k+2\ell\cr 0\cr} \cr
\noalign{\smallskip}
\lambda_{w_1} = (p-{\sqrt{5}\over6}q) - {q\over6}k\ell: \quad&
	w_1 = \pmatrix{3k\cr \sqrt{5}\,j-2\,i\ell\cr -7-\sqrt{5}\,k\ell\cr} \cr
\noalign{\smallskip}
\lambda_{w_2} = (p-{\sqrt{5}\over6}q) + {q\over6}k\ell: \quad&
	w_2 = \pmatrix{\sqrt{5}\,k+2l\cr 3j\cr -3\sqrt{5}-3\,k\ell\cr} \cr
\noalign{\smallskip}
  }\Eqno$$
However, we have been unable to find any decompositions of $\hat{\BB}$
involving these vectors.  It is intriguing that, for instance, $v_1$ is
orthogonal to both $u_1$ and $w_1$ (in the sense of \Ortho), but that $u_1$
and $w_1$ are not orthogonal.  In fact, we have shown using {\sl Mathematica}
that there is {\it no\/} eigenvector triple containing $w_1$ which is
orthogonal in the sense of \Ortho.  Unless $w_1$ is special in some as yet to
be determined sense, we are forced to conclude that neither \Ortho\ nor
\DecompGen\ are generally true for eigenvectors whose eigenvalues are not
real.  It is curious, however, that the sum of the squares (outer products) of
all {\it six} of these (normalized) vectors is indeed (twice) the identity!
We will consider possible implications of this fact below.

\Subsection{Example 3}

In all of the examples considered so far, the eigenvalues have been in the
complex subalgebra of $\OO$ determined by the associator $[a,b,c]$ (with $a$,
$b$, $c$ as in \Three).  We now give an example for which this is not the
case.

Consider
$$C = \pmatrix{~~p & iq & -q(j-i\ell-j\ell) \cr \noalign{\smallskip}
		-iq & p & ~q(1+k+l) \cr \noalign{\smallskip}
		q(j-i\ell-j\ell) & -q(1-k-l) & p \cr}
 \Eqno$$
which admits an eigenvector
$$v = \pmatrix{j\cr l\cr 0} \Eqno$$
with eigenvalue
$$\lambda_v = p + q \, lk \Eqno$$
However, the associator takes the form
$${[a,b,c] \over q^3} = [i,(j-i\ell-j\ell),(1+k+l)] = 2 (l-k) \Eqno$$

\Section{DISCUSSION
\Footnote{Much of this discussion also appeared in \Nonreal.}
}

As pointed out in \Eigen, the orthonormality relation \Ortho\ is equivalent to
assuming that 
$$vv^\dagger + ... + ww^\dagger=I \Eqno$$
If we define a matrix $U$ whose columns are just $v,...,w$, then this
statement is equivalent to
$$U U^\dagger = I \Eqno$$\label\UOrtho
Furthermore, the eigenvalue equation \Master\ can now be rewritten in the form
$$A U = U D \Eqno$$ \label\Matrix
where $D$ is a diagonal matrix whose entries are the eigenvalues.
Decompositions of the form \DecompSp\ now take the form
$$A = (UD) U^\dagger \Eqno$$
and multiplication of \Matrix\ on the right by $U^\dagger$ shows that
$$(AU) U^\dagger = (UD) U^\dagger = A = A (UU^\dagger) \Eqno$$
Thus, just as in \Eigen, decompositions of the form \DecompSp\ can be viewed
as the assertion of associativity
$$(AU) U^\dagger = A (UU^\dagger) \Eqno$$\label\UAssoc
We know of no way to express decompositions of the form \DecompGen\ in similar
language, which leads us to suspect that \DecompSp\ is more fundamental.

This is somewhat less surprising when one realizes that
$$\Big( (v \lambda) \, v^\dagger \Big) \, v = (v \lambda) \Eqno$$\label\ID
for {\it any\/} (normalized) vector $v$ and octonion $\lambda$, due to an
associator identity discussed in the Appendix of \Nonreal.  One can therefore
construct matrices with arbitrary octonionic (right) eigenvalues, although
such matrices will not be Hermitian.  We conjecture, however, that the correct
notion of orthogonality is
$$\Big( (v \lambda) \, v^\dagger \Big) \, w
  = 0 \Eqno$$\label\NewOrtho
not \Ortho, to which it of course reduces if the eigenvalues are real, in
which case it might be possible to decompose a matrix into non-Hermitian
pieces of the form \ID, and yet have the sum of these pieces be Hermitian.
This is precisely what happens in the $2\times2$ case \Nonreal, as well as in
Example~1.  In any case, it is intriguing that this notion of orthogonality
can be written as
$$\Big( (\AA v) \, v^\dagger \Big) \, w = 0 \Eqno$$
which explicitly involves $\AA$.

Putting these ideas together, it would be natural to conjecture that {\it all}
eigenvectors of a $3\times3$ octonionic Hermitian matrix come in families of
3, which form a decomposition in the sense that \UAssoc\ is satisfied, and
which are orthogonal in the sense of \NewOrtho.  However, this conjecture
appears to fail for Example~2, as our claim that there is no orthonormal (in
the sense \Ortho) eigenvector triple containing $w_1$ rules out a
decomposition (involving $w_1$) in either the form \DecompGen\ or the form
\DecompSp.  The reason for this is that the eigenvectors of $\hat\BB$ do not
depend on the parameters $p$ and $q$, which in turn means that the $p$ and $q$
parts of the decompositions can be treated separately.  If a decomposition of
either form were to exist, the terms involving $p$ would then imply \UOrtho,
contradicting our claim.  (This is the reason that we did not reduce to the
case $p=0$, $q=1$.)

Nonetheless, for all known decompositions \UOrtho\ does hold, and could be
used to further simplify \UAssoc.  One possible generalization would be for
\UOrtho\ to fail, but for a decomposition along the lines of \UAssoc\ to hold.
However, even \UAssoc\ appears to fail for Example~2.

There is, however, another intriguing possibility.  Example~2 suggests that
the eigenvectors of $3\times3$ octonionic Hermitian matrices may come in sets
of 6 (or more), rather than in sets of 3.  This would fit nicely with our
recent result with Okubo
\Ref{Tevian Dray, Corinne A. Manogue, and Susumu Okubo, in preparation.}
that, for real eigenvalues, it takes all 6 eigenvectors in order to decompose
an arbitrary vector into a linear combination of eigenvectors, despite the
fact that only 3 eigenvectors are needed to decompose the original matrix.
Further evidence for this point of view is provided by the fact that the most
general eigenvectors for the given eigenvalues of $\hat{\BB}$ have at most 2
free parameters, rather than the 4 degrees of freedom shown in \Eigen\ to
exist for real eigenvalues, or the 8 for complex matrices.

We therefore conjecture that, for any $3\times3$ octonionic Hermitian matrix,
\UAssoc\ should hold when suitably rewritten for a set consisting of $n$
eigenvectors, where $n$ presumably divides 24, the number of (real)
independent eigenvectors with real eigenvalues.  However, we have so far been
unable to write Example~2 in such a form.  Whether or in what form
orthogonality would hold in such a context is an interesting open question.

\vfill\eject
\References

\bye